\documentclass[12pt]{extarticle}

\usepackage[margin=1.5in]{geometry}  % set the margins to 1in on all sides
\usepackage{graphicx}              % to include figures
\usepackage{amsmath}               % great math stuff
\usepackage{amsfonts}              % for blackboard bold, etc
\usepackage{amsthm}                % better theorem environments
\usepackage{framed}
\usepackage{tikz}

\newtheorem{thm}{Theorem}[section]
\newtheorem{lem}[thm]{Lemma}

\begin{document}

\nocite{*}

\title{\bf Gaps Between Almost-Primes \\ and a Construction of\\ Almost-Ramanujan Graphs}

\author{\textsc{Adrian Dudek} \\ 
Mathematical Sciences Institute \\
The Australian National University \\ 
\texttt{adrian.dudek@anu.edu.au}}
\date{}

\maketitle

\begin{abstract}
For all $k \geq 3$, we show how one can explicitly construct an infinite family of $k$-regular graphs all of which have second largest eigenvalue satisfying the bound $O(k^{1/2})$. This resolves an open problem of Reingold, Vadhan and Wigderson. 
\end{abstract}

\section{Introduction}

In this note, we will consider a graph $X$ with vertex set $V$. By a family of graphs, we mean a sequence $\{X_m\}$ of graphs with $|V_m| \rightarrow \infty$ as $m \rightarrow \infty$. We list the eigenvalues of the adjacency matrix of a graph $X$ as
$$\lambda_1(X) \geq \lambda_2(X) \geq \cdots \geq \lambda_{|V|} (X).$$

Given a finite, connected, $k$-regular graph, it is a straightforward exercise to show that $\lambda_1 = k$. It is, however, of more interest to obtain upper bounds on $\lambda_2$; the following theorem of Alon and Boppana \cite{alon} gives a limit on this asymptotically. 

\begin{thm}
Let $\{X_m\}$ be a family of finite, connected, $k$-regular graphs. Then
$$\liminf_{m \rightarrow \infty} \lambda_2(X_m) \geq 2 \sqrt{k-1}.$$
\end{thm}

Thus, asymptotically, the spectral gap $\lambda_1 - \lambda_2$ is bounded above by $k-2 \sqrt{k-1}$. As usual, we define a \textit{Ramanujan graph} to be a finite, connected, $k$-regular graph with $\lambda_2 \leq 2 \sqrt{k-1}$. The central problem in the theory, given an integer $k \geq 3$, is to construct a $k$-regular family of Ramanujan graphs, more precisely a sequence $\{X_m\}$ of finite, connected, $k$-regular graphs with $\lambda_2(X_m) \leq 2 \sqrt{k-1}$ for all $m \geq 1$.

These would be the best possible expanders, for the spectral gap would be as large as is asymptotically possible. The centrepiece of the theory of expander graphs is that families of Ramanujan graphs have indeed been constructed for all $k = p^a+1$ where $p$ is a prime and $a$ is a positive integer (see Lubotzky, Phillips and Sarnak \cite{lps}, Chiu \cite{chiu}, Margulis \cite{margulis2} and Morgenstern \cite{morgenstern}). As such, the first case where a construction of a family of Ramanujan graphs is not known is $k=7$. 

In lieu of such a result, it is then of interest as to whether one can explicitly construct families whose spectral gap is of the form $k-O(k^{1/2})$. It is the purpose of this note to outline such a construction.

Reingold, Vadhan and Wigderson \cite{rvw} developed the so-called zig-zag product, and used this to explicitly construct $k$-regular families of graphs with second largest eigenvalues bounded above by $O(k^{2/3})$ for all $k \geq 3$. They then asked if one could improve this to $O(k^{1/2})$, which has the same order as that of the Ramanujan graphs. Ben-Aroya and Ta-Shma \cite{bt} came close and were able to explicitly construct $k$-regular families with
$$\lambda_2 \leq k^{1/2+O(1/\sqrt{\log k})}$$
In an earlier paper, the author \cite{dudek} showed how one could improve on this slightly assuming the truth of the Riemann hypothesis. One starts with the explicit construction of a family of $d$-regular Ramanujan graphs where $d$ is one more than a prime, and then uses the Cartesian graph product to produce a $(d+1)$-regular family of graphs with at least the same spectral gap. The reader should note that the new family is not necessarily Ramanujan, as we should expect an increase in the spectral gap with respect to an increase in $d$. We can iterate this procedure to furnish a family of $k$-regular graphs for any $k$. On the Riemann hypothesis, Cram\'{e}r \cite{cramer} proved that one can bound the difference between consecutive primes $p_n, p_{n+1}$ with
$$p_{n+1} - p_n = O(\sqrt{p_n} \log p_n),$$
and thus, one only needs to iterate so far. The method gave an explicit construction of $k$-regular families all of whose graphs satisfy
$$\lambda_2 \leq k^{1/2+O(\log \log k / \log k)}$$

In this note we resolve, without any condition, the open problem of Reingold, Vadhan and Wigderson \textit{viz.} the following theorem.

\begin{thm} \label{main}
For all sufficiently large integers $k$, one can explicitly construct a $k$-regular family $\{X_m\}$ of graphs such that
$$\lambda_2(X_m) \leq 4 \sqrt{k-1} + k^{101/232}$$
for all integers $m \geq 1$.
\end{thm}

It follows from the above theorem, that for all integers $k \geq 3$ we can explicitly construct a $k$-regular family of graphs with spectral gaps $k - O(k^{1/2})$.

\section{Proof of Theorem \ref{main}}

The proof is short, and is a novel twist of the author's earlier work, as one now starts with $k$-regular families where $k$ is a product of at most two primes and iterates.  Our construction thus relies on two important results. The first is the construction of almost-Ramanujan graphs given by Pizer \cite{pizer}; we state this in the form of the following lemma.

\begin{lem} \label{lemmapizer}
For $q \geq 3$, one can explicitly construct a family of $(q+1)$-regular graphs $\{X_m\}$ satisfying
$$\lambda_2 (X_m) \leq d(q+1) \sqrt{q}$$
for all $m \geq 1$, where $d(n)$ denotes the number of divisors of $n$.
\end{lem}

The well-known divisor bound states that 
$$d(n) \leq n^{O(1/\log \log n)},$$
and so one immediately has the explicit construction of a $(q+1)$-regular family $\{X_m\}$ of graphs such that
$$\lambda_2(X_m) \leq q^{1/2 + O(1/\log \log q)}$$
for all $m \geq 1$. This, however, falls short of the result of Ben-Aroya and Ta-Shma.

We will also need the following result of Wu \cite{wu} on almost-primes in short intervals.

\begin{lem} \label{lemmawu}
Let $x$ be sufficiently large. Then there exists an integer $n$ in the interval
$$(x-x^{101/232}, x]$$
such that the number of prime factors of $n$ counted with multiplicity does not exceed two.
\end{lem}

The proof proceeds as follows. Given two graphs $X$ and $Y$, the \textit{Cartesian product} $X \square Y$ is a natural way to obtain a new graph whose properties reflect those of the original graphs. We will not need to see the general definition of this, for we are using one specific instance of the product.

Consider $K_2$, the complete graph on two vertices, that is, the graph consisting of two vertices connected by a single edge. Given a finite $k$-regular graph $X$, we note that $X \square K_2$ is obtained simply by taking $X$ and a duplicate of $X$, and connecting each vertex in $X$ with its duplicate vertex. We will denote this particular Cartesian product by $X'$ and note that this is a $(k+1)$-regular graph with twice as many vertices as $X$.

The following result can be found in an earlier paper by the author \cite{dudek}.

\begin{thm} \label{up}
Let $X$ be a finite, connected $k$-regular graph and let $X' = X \square K_2$. Then
\begin{equation}
\lambda_2(X') \leq \lambda_2(X)+1.
\end{equation}
\end{thm}

Therefore, one may take a $k$-regular family $\{X_m\}$ of graphs and get a $(k+1)$-regular family $\{X_m'\}$ with the same upper bound on the second largest eigenvalue. Iterating this process will give families of all degrees.

Now, let $k$ be sufficiently large, so that we can define $q$ to be such that $q$ has at most two prime factors and
\begin{equation} \label{qbound}
k - k^{101/232} < q \leq k
\end{equation}
 where $p$ and $q$ are both prime numbers. By Lemma \ref{lemmapizer}, it follows that we can construct a family $\{X_m\}$ of $q$-regular graphs such that 
$$\lambda_2 (X_m) \leq 4 \sqrt{q-1}.$$
We now repeatedly take the Cartesian product of each graph in the family $\{X_m\}$ with $K_2$ so that we end up with a $k$-regular family $\{Y_m\}$. Note that by Theorem \ref{up} and (\ref{qbound}) we have 
$$\lambda_2(Y_m) \leq 4 \sqrt{q-1} + (k-q) < 4 \sqrt{k-1} + k^{101/232},$$
and this completes the proof of Theorem \ref{main}.

\clearpage

\bibliographystyle{plain}

\bibliography{biblio}

\end{document}